\documentclass{article}

\usepackage[pstarrows]{pict2e}
\usepackage{amsmath,amsfonts,a4wide,latexsym,amssymb,theorem,color,picture,graphicx}
\usepackage[pdftex, colorlinks, citecolor=blue]{hyperref}

\usepackage[T1]{fontenc}
\usepackage[utf8]{inputenc}
\usepackage{textcomp}
\usepackage{xcolor}
\usepackage{enumerate}
\usepackage{mathrsfs}
\usepackage{ulem}

\newcommand\A{{\mathscr{A}}}

\newcommand\K{{\mathbb K}}

\renewcommand\L{{\mathbb L}}

\newcommand\N{{\mathbb N}}

\newcommand\Q{{\mathbb Q}}
\newcommand\QQ{\bar{\mathbb Q}}

\newcommand\Z{{\mathbb Z}}

\newcommand\height{{\mathrm h}}
\newcommand\lgcd{{\mathrm {lgcd}}}
\newcommand\lcm{{\mathrm {lcm}}}

\newcommand{\ua}{{\underline{a}}}
\newcommand{\uw}{{\underline{w}}}
\newcommand{\ux}{{\underline{x}}}
\newcommand{\uy}{{\underline{y}}}
\newcommand{\uzero}{{\underline{0}}}

\newcommand{\tilf}{{\tilde f}}
\newcommand{\tilS}{{\widetilde S}}
\newcommand{\tilT}{{\widetilde T}}

\newcommand\qed{\hfill$\square$}

\theorembodyfont{\slshape}
\newtheorem{proposition}{Proposition}[section]
\newtheorem{theorem}[proposition]{Theorem}
\newtheorem{remark}[proposition]{Remark}
\newtheorem{lemma}[proposition]{Lemma}
\newtheorem{corollary}[proposition]{Corollary}

\title{The Skolem-Abouzaïd's theorem in the singular case}
\author{Boris Bartolome\\
Enteleia Tech\\
La Cour\\
31320 Aureville\\
France\\
Boris.Bartolome@enteleia-tech.com}

\date\today

1

\newcommand{\eps}{\varepsilon}

\begin{document}

\hfuzz4pt

\maketitle




\section{Introduction}

Let ${F(X,Y)\in\Q[X,Y]}$ be a $\Q$-irreducible polynomial.  In 1929 Skolem~\cite{Sk29} proved the following beautiful theorem:

\begin{theorem}[Skolem]
\label{tskol}
Assume that 
\begin{equation}
\label{e00}
F(0,0)=0.
\end{equation} 
Then for every non-zero integer~$d$, the equation ${F(X,Y)=0}$ has only finitely many solutions in integers ${(X,Y)\in \Z^2}$ with ${\gcd(X,Y)=d}$. 
\end{theorem}

In the same year, Siegel obtained his celebrated finiteness theorem for integral solutions of Diophantine equations: equation ${F(X,Y)=0}$ has finitely many solutions in integers unless the corresponding plane curve is of genus~$0$ and has at most~$2$ points at infinity. While Siegel's result is, certainly, deeper and more powerful than Theorem~\ref{tskol}, the latter has one important advantage. Siegel's theorem is known to be non-effective: it does not give any bound for the size of integral solutions. On the contrary, Skolem's method allows one to bound the solutions explicitly in terms of the coefficients of the polynomial~$F$ and the integer~$d$. Indeed, such a bound was obtained by Walsh~\cite{Walsh}; see also~\cite{Poulakis}.

In 2008, Abouzaid~\cite{Ab08} gave a far-going generalization of Skolem's theorem. He extended it in two directions.

First, he studied solutions not only in rational integers, but in arbitrary algebraic numbers. To accomplish this, he introduced the notion of \textsl{logarithmic gcd} of two algebraic numbers~$\alpha$ and~$\beta$,  which coincides with the logarithm of the usual gcd when ${\alpha, \beta\in  \Z}$. 

Second, he not only bounded the solution in terms of the logarithmic gcd, but obtained a sort of asymptotic relation between the heights of  the coordinates and their logarithmic gcd.  

Let us state Abouzaid's principal result (see~\cite[Theorem~1.3]{Ab08}). In the sequel we assume that ${F(X,Y) \in \bar\Q[X,Y]}$ is an absolutely irreducible polynomial, and use the notation
\begin{equation}
\label{emn}
m=\deg_XF, \quad n=\deg_YF, \quad M = \max\{m, n\}.
\end{equation}
We denote by $\height(\alpha)$ the absolute logarithmic height of ${\alpha\in \bar\Q}$ and by $\lgcd(\alpha,\beta)$ the logarithmic gcd of $\alpha,\beta\in \Q$. We also denote by $\height_p(F)$ the projective height of the polynomial~$F$. For all definitions, see Subsection~\ref{ssprem}.

\begin{theorem}[Abouza\"id]
Assume that $(0, 0)$ is a non-singular point of the plane curve $F(X,Y) = 0$. Let~$\eps$ satisfy ${0 < \eps <1}$. Then for any solution ${(\alpha,\beta)\in \bar\Q^2}$  of  ${F(X,Y) = 0}$, we have either 
$$
\max\{\height(\alpha), \height(\beta)\} \le  56M^8\eps^{-2}\height_p(F) + 420M^{10}\eps^{-2} \log(4M),
$$
 or
\begin{align*}
\max\{|\height(\alpha) - n\lgcd(\alpha,\beta)|, |\height(\beta) - m\lgcd(\alpha,\beta)|\} \le\ & \eps \max\{\height(\alpha), \height(\beta)\} + 742M^7\eps^{-1}\height_p(F)\\& + 5762M^9\eps^{-1} \log(2m + 2n).
\end{align*}
\end{theorem}

Informally speaking,
\begin{equation}
\label{easymp}
\frac{\height(\alpha)}n\sim\frac{\height(\beta)}m\sim\lgcd(\alpha,\beta)
\end{equation}
as ${\max\{\height(\alpha), \height(\beta)\}\to \infty}$.

Unfortunately, Abouzaid's assumption is slightly more restrictive than Skolem's~\eqref{e00}: he assumes not only that the point $(0,0)$ belongs to the plane curve ${F(X,Y)=0}$, but also that $(0,0)$ is a  non-singular point on this curve. 

Denote by~$r$ the ``order of vanishing'' of $F(X,Y)$ at the point $(0,0)$:
\begin{equation}
\label{eordvan}
r=\min\left\{i+j: \frac{\partial^{i+j}F}{\partial^iX\partial^jY}(0,0)\ne0\right\}. 
\end{equation}
Clearly, ${r>0}$ if and only if  ${F(0,0)=0}$ and  ${r=1}$ if and only $(0,0)$ is a non-singular point of the plane curve ${F(X,Y)=0}$. 

We can now state our principal result. 

\begin{theorem}
\label{MainResult}
Let ${F(X,Y)\in \bar\Q[X,Y]}$ be an absolutely irreducible polynomial satisfying $F(0,0)=0$. 
Let~$\eps$ satisfy ${0 < \eps <1}$.
Then, for any $\alpha, \beta \in \QQ$ such that $F(\alpha, \beta)=0$, we have either:
\begin{equation*}
\height (\alpha) \leq 200\eps^{-2}mn^6(\height_p(F)+5)
\end{equation*}
or
\begin{equation*}
\left| \frac{\lgcd (\alpha, \beta)}{r} - \frac{\height(\alpha)}{n} \right| \leq
\frac{1}{r} \left(   \eps\height(\alpha)+4000\eps^{-1} n^4(\height_p(F)+\log(mn) + 1) +30n^2m(\height_p (F) + \log (nm))  \right).
\end{equation*}
\end{theorem}

By symmetry, the same kind of bound holds true for the difference ${ \frac{\lgcd (\alpha, \beta)}{r} - \frac{\height(\beta)}{m} }$. Informally speaking,
\begin{equation}
\label{einformal}
\frac{\height(\alpha)}n\sim\frac{\height(\beta)}m\sim\frac{\lgcd(\alpha,\beta)}r
\end{equation}
as ${\max\{\height(\alpha), \height(\beta)\}\to \infty}$. 

Validity of~\eqref{einformal} was conjectured by Abouzaid, see the end of Section~1 in~\cite{Ab08}\footnote{Abouzaid's definition of~$r$ looks different, but it can be easily shown that it is equivalent to ours.}. 
The referee  pointed  us to an unpublished work of Habegger~\cite{Hab} from 2007, where he confirms Abouzaid's conjecture; moreover, his bounds are sharper than ours. We would like to remark that Habegger's method is quite different and uses his sharp quantitative version of the quasi-equivalence of heights. On the contrary, our    paper follows closely the methods of~\cite{Ab08} wherever possible; in particular, like in~\cite{Ab08}, our main tool is Puiseux expansions.

\paragraph{Plan of the article}
Section~\ref{sheights} and~\ref{sseries} are preliminary: we compile therein some definitions and results from different sources, which will be used in the article. In Section~\ref{sspr} we establish the ``Main Lemma'', which is the heart of the proof of Theorem~\ref{MainResult}. In Section~\ref{sproof} we complete the proof of Theorem~\ref{MainResult} using the ``Main Lemma''.

\paragraph{Acknowledgments}

I am grateful to Yuri Bilu for having pointed my attention to this problem and for an emulating exchange on this topic. I am also thankful to the referee for her/his helpful suggestions and for pointing out the unpublished result from Philip Habbeger.

\section{Heights}
\label{sheights}
In this section we recall definitions and collect various results about absolute values and heights. 

We normalize the absolute values on number fields so that they extend standard absolute values on~$\Q$: if ${v\mid p}$ (non-Archimedean) then ${|p|_v=p^{-1}}$ and if ${v\mid \infty}$ (Archimedean) then ${|2015|_v=2015}$. 

\subsection{Heights and lgcd of algebraic numbers}
\label{ssprem}

Let~$\K$ be a number field, $d = [\K:\Q]$ and $d_v = [\K_v:\Q_v]$. The \textsl{height} of an algebraic number~$\alpha \in \K$ is defined as
$$
\height(\alpha) = \frac1{d}\sum_{v\in M_\K}d_v\log^+|\alpha|_v.
$$
where $M_\K$ is the set of places (normalized absolute values) of the number field~$\K$ and $\log^+=\max\{\log,0\}$. It is well-known that the height does not depend on the particular choice of~$\K$, but only on the number~$\alpha$ itself. It is equally well-known that 
${\height(\alpha)=\height(\alpha^{-1})}$, so that 
$$
\height(\alpha) = \frac1{d}\sum_{v\in M_\K}-d_v\log^-|\alpha|_v = \sum_{v\in M_\K}\height_v(\alpha),
$$
where  ${\log^-=\min\{\log,0\}}$ and 
$$
\height_v(\alpha)=-\frac{d_v}{d}\log^-|\alpha|_v. 
$$
The quantities $\height_v(\alpha)$ can be viewed as ``local heights''. Clearly, ${\height_v(\alpha)\ge 0}$ for any~$v$ and~$\alpha$. 

We define the \textsl{logarithmic gcd} of two algebraic numbers~$\alpha$ and~$\beta$, not both~$0$, as 
$$
\lgcd (\alpha,\beta) =\sum_{v\in M_K}\min\{\height_v(\alpha), \height_v(\beta)\},
$$
where~$\K$ is a number field containing both~$\alpha$ and~$\beta$. It again depends only~$\alpha$ and~$\beta$, not on~$\K$. A simple verification shows that for ${\alpha, \beta \in \Z}$ we have 
${\lgcd(\alpha,\beta)=\log\gcd(\alpha,\beta)}$.

Now let~$\K$ be a number field and~$S$ be a set of places of~$\K$. We define the \textsl{$S$-height} by 
$$
\height_S(\alpha) = \sum_{v\in S}\height_v(\alpha).
$$
Similarly we define $\lgcd_S$. We shall frequently use the inequality $\lgcd_S (\alpha, \beta) \leq \height_S (\alpha) \leq \height(\alpha)$ without special reference.

\subsection{Affine and projective heights of polynomials}

We define the projective and the affine height of a vector $\underline{a}=(a_1, \ldots, a_m) \in \QQ^m$ with algebraic entries, by
\begin{align*}
\height _p (\underline{a}) &= \frac{1}{d}\sum_{v \in M_{\K}} d_v \log \max_{1 \leq k \leq m} |a_k|_v\qquad (\underline{a} \neq \underline{0}),\\
\height _a (\underline{a}) &= \frac{1}{d}\sum_{v \in M_{\K}} d_v \log^+ \max_{1 \leq k \leq m} |a_k|_v.
\end{align*}
Here, $\K$ is a number field containing $a_1, \ldots, a_m$, and $d$, $d_v$ are defined as in the previous subsection. We notice that the height of an algebraic number defined in the previous subsection corresponds to the affine height of a one-dimensional vector.

We define the projective and affine height of a polynomial as the corresponding heights of the vector of its non-zero coefficients. If $F$ is a non-zero polynomial, then, for ${\alpha \in \QQ^{*}}$ we have ${\height_p (\alpha F) = \height_p(F)}$. Also, ${\height_p ( F) \leq \height_a(F)}$, with  ${\height_p ( F) = \height_a(F)}$ if $F$ has a coefficient equal to $1$.

In \cite[Lemma 4]{Schmidt2}, Schmidt proves the following lemma:

\begin{lemma}
\label{lSchmidt2}
Let $F(X,Y) \in \QQ[X,Y]$ be a polynomial with algebraic coefficients, such that $m = \deg_X F$ and $n = \deg_Y F$. Let $R_F(X) = \text{Res}_Y (F, F_Y^{'})$ be the resultant of F and its derivative polynomial with respect to Y. Then:
\begin{equation}
\height_p(R_F) \leq (2n-1)\height_p (F) + (2n-1) \log ((m+1)(n+1)\sqrt{n}).
\end{equation}
\end{lemma}

It is well-known that the height of a root of a polynomial is bounded in terms of the height of the polynomial itself. The following lemma can be found in 
\cite[Proposition 3.6]{Yuri2}:

\begin{lemma}
\label{Yuri2}
Let $F(X)$ be a polynomial of degree $m$ with algebraic coefficients. Let $\alpha$ be a root of $F$. Then, $\height (\alpha) \leq  \height_p(F) + \log 2$
\end{lemma}

We want to generalize this to a system of two algebraic equations in two variables.

\begin{lemma}
\label{lsys}
Let $F_1(X,Y)$ and $F_2(X,Y)$ be polynomials with algebraic coefficients, having no common factor.  Put:
\begin{equation*}
m_i=\deg_X F_i,\ n_i = \deg_Y F_i \qquad (i=1,2).
\end{equation*}
Let $\alpha,\beta$ be algebraic numbers satisfying ${F_1(\alpha,\beta)=F_2(\alpha,\beta)=0}$. 
Then
\begin{equation*}
\height (\alpha) \leq n_1 \height_p (F_2) + n_2 \height_p (F_1) + (m_1 n_2 + m_2 n_1) + (n_1 + n_2) \log (n_1 + n_2)+\log2.
\end{equation*}
\end{lemma}

\paragraph{Proof}
Since~$F_1$ and~$F_2$ have no common factor, their $Y$-divisor $R(X)$ is a non-zero polynomial, and ${R(\alpha)=0}$. \cite[Proposition 2.4]{Ab08} gives the estimate 
\begin{equation*}
\height_p (R) \leq n_1 \height_p (F_2) + n_2 \height_p (F_1) + (m_1 n_2 + m_2 n_1) + (n_1 + n_2) \log (n_1 + n_2).
\end{equation*}
Combining this with Lemma~\ref{Yuri2}, the result follows. \qed

\bigskip

We will also use \cite[Proposition 2.5]{Ab08}:
\begin{lemma}
\label{lmaxG}
Let $F(X,Y) \in \QQ[X,Y]$ be a polynomial with $m = \deg_X F$ and $n = \deg_Y F$ and let $\alpha, \beta$ be two algebraic numbers. Then
\begin{enumerate}
\item 
\label{ifab}
We have $\height (F(\alpha, \beta)) \leq \height_a (F) + m\height(\alpha) +n\height(\beta) + \log((m+1)(n+1))$.
\item 
\label{iba}
If $F(\alpha, \beta) = 0$ with $F(\alpha, Y)$ not vanishing identically, then:
\begin{equation*}
\height(\beta) \leq \height_p(F) + m\height(\alpha) +n +\log (m+1).
\end{equation*}
\end{enumerate}
\end{lemma}

\subsection{Coefficients versus roots}

In this subsection we establish some simple relations between coefficients and roots of a polynomial over a field with absolute value, needed in the proof of our main result. It will be convenient to use the notion of \textsl{$v$-Mahler measure} of a polynomial. 

Let~$\K$ be a field with absolute value~$v$ and ${f(X)\in \K[X]}$ a  polynomial of degree~$n$. Let  ${\beta_1, \ldots, \beta_n\in \bar \K}$ be the roots of $f$:
$$
f(X)=a_nX^n+a_{n-1}X^{n-1}+\ldots+a_0=a_n(X-\beta_1)\ldots(X-\beta_n).
$$
Define the $v$-Mahler measure of~$f$ by
$$
M_v(f)=|a_n|_v\prod_{i=1}^n\max\{1,|\beta_i|_v\}, 
$$
where we extend~$v$ somehow to~$\bar \K$. (Clearly, $M_v(f)$ does not depend on the particular extension of~$v$.)
It is well-known that ${|f|_v=M_v(f)}$ for non-archimedean~$v$ (``Gauss lemma'') and ${M_v(f)\le (n+1)|f|_v}$ for archimedean~$v$ (Mahler).

\begin{lemma}
\label{lr+1}
Let ${\beta_1, \ldots, \beta_{\ell+1}}$ be ${\ell+1}$ distinct roots of $f(X)$, where ${0\le \ell\le n-1}$. Then
$$
\max\{|\beta_1|_v, \ldots, |\beta_{\ell+1}|_v\}\ge c_v(n)\frac{|a_\ell|_v}{|f|_v},
$$
where ${c_v(n)=1}$ for non-archimedean~$v$ and ${c_v(n)=(n+1)^{-1}2^{-n}}$ for archimedean~$v$.
\end{lemma}

\paragraph{Proof}
We have
\begin{equation}
\label{eanu}
a_\ell=\pm a_n\sum_{1\le i_1<\ldots<i_{n-\ell}\le n}\beta_{i_1}\ldots\beta_{i_{n-\ell}},
\end{equation}
where ${\beta_1, \ldots, \beta_n}$ are all roots of $f(X)$ in~$\bar\K$ counted with multiplicities. 
Observe that each term in the sum above contains one of the roots ${\beta_1, \ldots, \beta_{\ell+1}}$, and the product of the other roots together with~$a_n$ is $v$-bounded by $M_v(f)$. 
Hence, denoting ${\mu= \max\{|\beta_1|_v, \ldots, |\beta_{\ell+1}|_v\}}$, we obtain 
${|a_\ell|_v\le \mu M_v(f)}$ in the non-archimedean case and ${|a_\ell|_v\le \binom n\ell\mu M_v(f)}$ in the archimedean case. Since ${\binom n\ell\le 2^n}$, the result follows. \qed

\subsection{Siegel's ``Absolute'' Lemma}

In this section we give a version of the Absolute Siegel's Lemma due to David and Philippon~[3], adapted for our purposes.

We start from a slightly modified definition of the projective height of a non-zero vector  ${\ua=(a_1, \ldots, a_n)\in\bar\Q^n}$. As before, we fix a number field~$\K$ containing ${a_1,\ldots,a_n}$ and set ${d=[\K:\Q]}$, \ ${d_v=[\K_v:\Q_v]}$ for ${v\in M_\K}$.

Now we define 
$$
\height_s(\ua)= \sum_{v\in M_\K}\frac{d_v}d \log \|\ua\|_v,
$$
where 
$$
\|\ua\|_v=
\begin{cases}
\max\{ |a_1|_v, \ldots, |a_n|_v\}, & v<\infty,\\
(|a_1|_v^2+\ldots+|a_n|_v^2)^{1/2}, & v\mid\infty.
\end{cases}
$$
This definition is the same as for $\height_p(\ua)$, except that for the archimedean places the sup-norm is replaced by the euclidean norm. 
We have clearly ${\height_s(\lambda\ua)=\height_s(\ua)}$ for ${\lambda\in \bar\Q^\times}$, and 
\begin{equation}
\label{ehphs}
\height_p(\ua)\le \height_s(\ua)\le \height_p(\ua)+ \frac12\log n. 
\end{equation}

Now let us define the height of a  linear subspace of $\bar\Q^n$.  If~$W$ is a $1$-dimensional subspace of $\bar\Q^n$ then we set 
$$
\height_s(W):=\height_s(\uw),
$$
where~$\uw$ is an arbitrary non-zero vector from~$W$. Clearly, $\height_s(W)$ does not depend on the particular choice of the vector~$\uw$. 

To extend this to subspaces of arbitrary dimension, we use Grassmann spaces. Recall that the $m$th Grassmann space 
$\wedge^m\bar\Q^n$ is of dimension $\binom nm$, and has a standard basis consisting of the vectors 
$$
e_{i_1}\wedge\ldots \wedge e_{i_m}, \qquad (1\le i_1<\ldots<i_m\le n), 
$$
where ${e_1, \ldots, e_n}$ is the standard basis of $\bar\Q^n$.  If~$W$ is an $m$-dimensional subspace of $\bar\Q^n$ then $\wedge^mW$ is a $1$-dimensional subspace of   $\wedge^m\bar\Q^n$, and we simply define
$$
\height_s(W):=\height_s(\wedge^mW). 
$$
Finally, we set ${\height_s(W)=0}$ for the zero subspace ${W=\{\uzero\}}$. 

To make this more explicit, pick a basis ${\uw_1,\ldots, \uw_m}$ of~$W$. Then $\wedge^mW$ is generated by ${\uw_1\wedge\ldots\wedge \uw_m}$, and we have 
\begin{equation}
\label{ehwedge}
\height_s(W)=\height_s(\uw_1\wedge\ldots\wedge \uw_m). 
\end{equation}
This allows one to estimate the height of a subspace generated by a finite set of vectors in terms of heights of generators.

\begin{proposition}
\label{elever}
Let~$W$ be a subspace of~$\bar\Q^n$ generated by vectors ${\uw_1,\ldots, \uw_m\in \bar\Q^n}$. Then
$$
\height_s(W)\le \height_s(\uw_1)+\ldots+\height_s(\uw_m).
$$
\end{proposition}

\paragraph{Proof}
Selecting among  ${\uw_1,\ldots, \uw_m}$ a maximal linearly independent subset, we may assume that ${\uw_1,\ldots, \uw_m}$ is a basis of~$W$. Then we have~\eqref{ehwedge}. It remains to observe that for any place~$v$ we have 
$$
\|\uw_1\wedge\ldots\wedge \uw_m\|_v\le \|\uw_1\|_v\ldots\|\uw_m\|_v. 
$$
For non-archimedean~$v$ this is obvious, and for archimedean~$v$ this is the classical Hadamard's inequality. \qed

\bigskip

We denote by $(\ux\cdot\uy)$ the standard inner product on~$\bar\Q^n$:
$$
(\ux\cdot\uy)=x_1y_1+\ldots+x_ny_n. 
$$ 
Let~$W^\perp$ denote the orthogonal complement to~$W$ with respect to this  product. It is well-known that the coordinates of $\wedge^mW$ (where ${m=\dim W}$) in the standard basis of $\wedge^m\bar\Q^n$ are the same (up to a scalar multiple) as the coordinates of $\wedge^{n-m}W^\perp$  in the standard basis of $\wedge^{n-m}\bar\Q^n$. In particular,
\begin{equation}
\label{eh=h}
\height_s(W)=\height_s(W^\perp).
\end{equation}
We use this to estimate the height of the subpace defined by a system of linear equations.

\begin{proposition}
Let $L_1, \ldots, L_m$ be non-zero linear forms on $\bar\Q^n$, and let~$W$ be the subspace of $\bar\Q^n$ defined by ${L_1(\ux)=\ldots=L_m(\ux)=0}$. Then
\begin{equation}
\label{ehwllp}
\height_s(W)\le \height_p(L_1)+\ldots+\height_p(L_m)+\frac m2\log n.
\end{equation}
\end{proposition}

\paragraph{Proof}
Let ${\ua_1, \ldots, \ua_m}$ be vectors in $\bar\Q^n$ such that ${L_i(\ux)=(\ux\cdot\ua_i)}$. Then 
\begin{equation}
\label{eal}
\height_p(L_i)=\height_p(\ua_i) \qquad (i=1, \ldots, m). 
\end{equation}
The space~$W^\perp$ is generated by ${\ua_1, \ldots, \ua_m}$. Applying to it Proposition~\ref{elever} and using~\eqref{ehphs}, we obtain 
$$
\height_s(W^\perp)\le \height_s(\ua_1)+\ldots+\height_s(\ua_m) \le \height_p(\ua_1)+\ldots+\height_p(\ua_m) +\frac m2\log n.
$$
Together with~\eqref{eh=h} and~\eqref{eal}, this gives~\eqref{ehwllp}. \qed

\begin{remark}
It is not difficult to slightly refine~\eqref{ehwllp}, replacing $\log n$ by $\log m$ in the right-hand side, but this would not lead to any substantial improvement of our results. 
\end{remark}

In [3, Lemma~4.7] the following version of ``absolute Siegel's lemma'' is given.

\begin{proposition}
\label{pdp}
Let~$W$ be an $\ell$-dimensional subspace of $\bar\Q^n$ and ${\eps>0}$. Then, there is a non-zero vector ${\ux\in W}$, satisfying:
\begin{equation*}
\height_p(\ux) \le \frac{\height_s(W)}{\ell}+ \frac1{2 \ell}\sum_{i=1}^{\ell-1}\sum_{k=1}^i\frac{1}{k}+\eps.
\end{equation*}
\end{proposition}

\begin{corollary}
\label{cdp}
Let $L_1, \ldots, L_m$ be non-zero linear forms in~$n$ variables with algebraic coefficients. Then, there exists a non-zero vector ${\ux\in \bar\Q^n}$ such that ${L_1(\ux)=\ldots=L_m(\ux)=0}$ and 
\begin{equation}
\label{eberny}
\height_p(\ux) \le \frac{1}{n-m}\left(\height_p(L_1)+\ldots+\height_p(L_m)\right)+ \frac12\frac{n}{n-m}\log n.
\end{equation}
\end{corollary}

\paragraph{Proof}
We apply Proposition~\ref{pdp} with~$W$ the subspace defined by ${L_1(\ux)=\ldots=L_m(\ux)=0}$. Denoting ${\ell=\dim W}$, we have clearly ${n-m\le r\le n}$ and
$$
\frac1{2\ell}\sum_{i=1}^{\ell-1}\sum_{k=1}^i\frac{1}{k}<\frac12\log \ell \le \frac12\log n. 
$$
Hence there exists a non-zero ${\ux\in W}$ satisfying 
$$
\height_p(\ux) \le \frac{1}{n-m}\height_s(W)+ \frac12\log n.
$$
Using~\eqref{ehwllp}, we find
$$
\height_p(\ux) \le \frac{1}{n-m}\left(\height_p(L_1)+\ldots+\height_p(L_m)\right)+\frac12\frac{m}{n-m}\log n+\frac12\log n,
$$
which is~\eqref{eberny}. \qed

\section{Power series}

\label{sseries}

In this section we recall various results about power series, used in our proof.

\subsection{Puiseux Expansions}
Let~$\K$ be a field of characteristic~$0$, and $\K((x))$ the field of formal power series over~$\K$. It is well-known that an  extension of $\K((x))$ of degree~$n$ is a subfield of a field of the form $\L((x^{1/e}))$, where~$e$ is a positive integer (the ramification index),~$\L$ is a finite extension of~$\K$, and 
$$
[\L:\K], e\le n.
$$ 
This fact (quoted sometimes as the ``Theorem of Puiseux'') has the following consequence:  if we fix an algebraic closure~$\bar \K$ of~$\K$, then the algebraic closure of $\K((x))$ can be given by
$$
\overline{\K((x))}=\bigcup_{e=1}^\infty\bigcup_{\genfrac{}{}{0pt}{}{\K\subset \L\subset \bar \K}{[\L:\K]<\infty}}\L((x^{1/e})),
$$
where the interior union is over all subfields~$\L$ of~$\bar \K$ finite over~$\K$. 

Another immediate consequence  of  the ``Theorem of Puiseux'' is the following statement:

\begin{proposition}
\label{ppe}
Let
$$
F(X,Y)=f_n(X)Y^n+\cdots+f_0(X)\in \K[X,Y]
$$ 
be a polynomial of $Y$-degree~$n$. Then there exists a finite extension~$\L$ of~$\K$, positive integers ${e_1, \ldots, e_n}$, all not exceeding~$n$, and series ${y_i\in \L((x^{1/e_i}))}$ such that 
\begin{equation}
\label{edecompf}
F(x,Y)=f_n(x)(Y-y_1)\cdots(Y-y_n).
\end{equation}
\end{proposition}

Write the series ${y_1,\ldots,y_n}$ as
$$
y_i=\sum_{k=\kappa_i}^\infty a_{ik}x^{k/e_i}
$$
with ${a_{i\kappa_i}\ne 0}$. It is well-known and easy to show that 
$$
|\kappa_i|\le \deg_XF  \qquad (i=1,\ldots, n). 
$$
This inequality will be used throughout the article without special notice. 

We want to link the numbers~$e_i$ and~$\kappa_i$ with the ``order of vanishing'' at $(0,0)$, introduced in~\eqref{eordvan}.

\begin{proposition}
\label{prke}
Let ${F(X,Y)\in \K[X,Y]}$ and ${y_1, \ldots, y_n}$ be as above, and assume that $F(0,Y)$ is not identically~$0$. Then the quantity~$r$, introduced in~\eqref{eordvan}, satisfies
\begin{equation}
\label{erke}
r=\sum_{\kappa_i>0}\min\{1, \kappa_i/e_i\},
\end{equation}
where the sum extends only to those~$i$ for which ${\kappa_i>0}$. 
\end{proposition}

\paragraph{Proof}
We denote by~$\nu_x$ the standard additive valuation on $\K((x))$, normalized to have ${\nu_x(x)=1}$. This~$\nu_x$ extends in a unique way to the algebraic closure $\overline{\K((x))}$; precisely, for  
$$
y(x)=\sum_{k=\kappa}^\infty a_kx^{k/e}\in \overline{\K((x))} \qquad (a_\kappa\ne0)
$$
we have ${\nu_x(y)=\kappa/e}$.
Furthermore, for  
$$
G(x,Y)=g_s(x)Y^s+\cdots+g_0(x)\in \overline{\K((x))}[Y]
$$
we set ${\nu_x(G)=\min \{\nu_x(g_0), \ldots, \nu_x(g_s)\}}$. 
Gauss' lemma asserts that for ${G_1,G_2 \in \overline{\K((x))}[Y]}$, we have ${\nu_x(G_1G_2)= \nu_x(G_1)+\nu_x(G_2)}$. 

Since $F(0,Y)$ is not identically~$0$, we have ${\nu_x(F(x,Y))=0}$. Applying Gauss' lemma to~\eqref{edecompf}, we obtain 
$$
\nu_x(f_0(x))+\sum\min\{0,\kappa_i/e_i\}=0.
$$
Hence, setting ${\tilf_0=x^{-\nu_x(f_0(x))}f_0(x)}$, we may re-write~\eqref{edecompf} as
\begin{equation}
\label{edecf}
F(x,Y)= \prod_{\kappa_i>0}(Y-y_i)\cdot \tilf_0(x)\prod_{\kappa_i\le 0}(x^{-\kappa_i/e_i}Y-x^{-\kappa_i/e_i}y_i). 
\end{equation}
Now set ${G(x,Y)=F(x,xY)}$. Then clearly ${r=\nu_x(G)}$. Applying Gauss' Lemma to the decomposition
$$
G(x,Y)= \prod_{\kappa_i>0}(xY-y_i)\cdot \tilf_0(x)\prod_{\kappa_i\le 0}(x^{1-\kappa_i/e_i}Y-x^{-\kappa_i/e_i}y_i), 
$$
we obtain~\eqref{erke}. \qed

\bigskip

Here is one more useful property. 

\begin{proposition}
\label{pfs}
In the set-up of Proposition~\ref{prke}, assume that ${\kappa_i>0}$ for exactly~$\ell$ indexes ${i\in \{1,\ldots, n\}}$. Then 
${f_k(0)=0}$ for ${k<\ell}$, but ${f_\ell(0)\ne 0}$. 
\end{proposition}

\paragraph{Proof}
Re-write~\eqref{edecf} as
$$
F(x,Y)= \prod_{\kappa_i>0}(Y-y_i)\prod_{\kappa_i=0}(Y-y_i)\cdot \tilf_0(x)\prod_{\kappa_i< 0}(x^{-\kappa_i/e_i}Y-x^{-\kappa_i/e_i}y_i). 
$$
Substituting ${x=0}$, every factor in the first product becomes~$Y$, every factor in the second product becomes ${Y-a_{i0}}$, with ${a_{i0}\ne 0}$, and every factor in the third product (including $\tilf_0(0)$) becomes constant. Whence the result. \qed. 

\subsection{Eisenstein's theorem}

In this subsection, we recall the quantitative Eisentsein's theorem due to work from Dwork, Robba, Schmidt and Van der Poorten \cite{DR79,DP92,Schmidt2},  as given in \cite{Yuri2}. It will be convenient to use the notion of $M_{\K}$-divisor.

An $M_{\K}$-divisor is an infinite vector $(A_v)_{v \in M_{\K}}$ of positive real numbers, each $A_v$ being associated to one $v \in M_{\K}$, such that for all but finitely many $v \in M_{\K}$ we have $A_v = 1$. An $M_{\K}$-divisor is effective if for all  $v \in M_{\K}, A_v \geq 1$. 

We define the \textit{height} of an $M_{\K}$-divisor $\A = (A_v)_{v \in M_{\K}}$ as

\begin{equation}
\height (\A) = \sum_{v \in M_{\K}} \frac{d_v}{d}\log A_v.
\end{equation}

The following version of Eisenstein's theorem is from \cite[Theorem 7.5]{Yuri2}. 

\begin{theorem}
Let $F(X,Y)$ be a separable polynomial of degrees ${m=\deg_X F}$ and ${n=\deg_Y F}$. Further, let $y(x)=\sum_{k = \kappa}^{\infty} a_{k} x^{k/e} \in \K [[x^{1/e}]]$ be a power series satisfying $F(x,y(x))=0$. (Here we do not assume that ${a_\kappa\ne 0}$.) Then there exists an effective $M_{\K}$-divisor $\A = \left( A_v \right)_{v \in M_{\K}}$ such that:
\begin{equation*}
|a_k|_v \leq \max \{1, |a_{e \lfloor \kappa /e \rfloor }|_v  \}A_v^{k/e - \lfloor \kappa/e \rfloor},
\end{equation*}
for any $v \in M_{\K}$ and any $k \geq \kappa$, and such that ${\height (\A) \leq 4n\height_p (F) + 3n \log (nm) + 10en}$.
\end{theorem}

Applying this theorem to the series of the form ${a_1x^{1/e}+a_2x^{2/e}+\ldots}$  (that is,  with ${a_k = 0}$ for ${k \leq 0}$) and setting $\kappa=0$, we obtain that:

\begin{corollary}
\label{cBB}
Let $F(X,Y)$ be a separable polynomial of degrees ${m=\deg_X F}$ and ${n=\deg_Y F}$. Further, let $y(x)=\sum_{k = 1}^{\infty} a_{k} x^{k/e} \in \K [[x^{1/e}]]$ be a power series satisfying $F(x,y(x))=0$.
Then, there exists an effective $M_{\K}$-divisor $\A = \left( A_v \right)_{v \in M_{\K}}$ such that:
\begin{equation}
\label{eeis1}
|a_k|_v \leq A_v^{k/e } \qquad(v\in M_\K, \quad k=1, 2, \ldots),
\end{equation}
and such that 
\begin{equation}
\label{eeis2}
\height (\A) \leq 4n\height_p (F) + 3n \log (nm) + 10en.
\end{equation}
\end{corollary}

The following lemma  is a slightly modified version of Proposition~2.7 from~\cite{Ab08}:

\begin{lemma}
\label{lPS}
Let $\K$ be a number field and let $y(x) = \sum_{k=1}^{\infty} a_k x^{k/e}$ be a series with coefficients in $\K$. Assume further that there exists an effective $M_{\K}$-divisor ${\A = \left( A_v \right)_{v \in M_{\K}}}$, such that for all  $k\geq 1$ we have $|a_k|_v \leq A_v^{k/e}$. For  $\ell \in \N$ write $y(x)^\ell = \sum_{k=1}^{\infty} a_k^{(\ell)} c^{k/e}$. Then, for any  ${v\in M_K}$ and for all  $k \geq 1$ we have:
\begin{equation}
|a_k^{(\ell)}|_v \leq 
\left\{   
\begin{array}{lcl}
2^{\ell+k} A_v^{k/e}, & \text{    } & \text{if } v | \infty,\\
A_v^{k/e}, & \text{          } & \text{if } v < \infty.
\end{array}
\right.
\end{equation}
\end{lemma}

In~\cite{Ab08}, a slightly sharper estimate, with   $\binom{\ell + k - 1}{k}$ instead of $2^{\ell+k}$ is given. 


\section{The ``Main Lemma''}
\label{sspr}

In this section we prove an auxiliary statement which is crucial for the proof of Theorem~\ref{MainResult}.  It can be viewed as a version of the famous Theorem of Sprindzhuk, see \cite{Bomb,BM06}. In fact, our argument is an adaptation of that from~\cite{BM06}. We follow \cite[Sections 3.1--3.3]{Ab08} with some changes. 

\subsection{Statement of the Main Lemma}
\label{ssml}
In this section~$\K$ is a number field, ${F(X,Y)\in \K[X,Y]}$ an absolutely irreducible  polynomial of degrees ${m=\deg_XF}$ and ${n=\deg_YF}$, and ${\alpha, \beta\in \K^\times}$ satisfy ${F(\alpha,\beta)=0}$. Furthermore, everywhere in this section except Subsection~\ref{ssrml}
$$
y(x)=\sum_{k=1}^\infty a_kx^k\in \K[[x]]
$$
is a power series satisfying ${F(x,y(x))=0}$; in particular, ${F(0,0)=0}$. 

We consider the following finite subset of $M_K$:
\begin{equation*}
T=\{v\in M_\K: \text{${|\alpha|_v<1}$ and $y(x)$ converges $v$-adically to~$\beta$ at ${x=\alpha}$}\}.
\end{equation*}

\begin{lemma}[``Main Lemma'']
\label{lml}
Let~$\eps$ satisfy ${0<\eps\le 1}$. Then  we have either
\begin{equation}
\label{eml1}
\height(\alpha)\le 200\eps^{-2}mn^4(\height_p(F)+5),
\end{equation}
or
\begin{equation}
\label{eml2}
\left|\frac{\height(\alpha)}n-\height_T(\alpha)\right|\le \eps n\height(\alpha)+200\eps^{-1} n^2(\height_p(F)+\log(mn) + 10). 
\end{equation}
\end{lemma}

\subsection{Preparations}

The proof of the ``Main Lemma'' requires some preparation. First of all, recall  that, according to Eisenstein's Theorem as given in Corollary~\ref{cBB}, there exists an effective $M_\K$-divisor ${\A = \left( A_v \right)_{v \in M_{\K}}}$ such that both~\eqref{eeis1} and~\eqref{eeis2} hold with ${e=1}$:
\begin{align*}
|a_k|_v &\leq A_v^{k} \qquad(v\in M_\K, \quad k=1, 2, \ldots),\\
\height (\A) &\leq 4n\height_p (F) + 3n \log (nm) + 10n.
\end{align*}
We fix this~$\A$ until the end of the section. 

Next, we need to construct an ``auxiliary polynomial''. 

\begin{proposition}[Auxiliary polynomial]
\label{pAuxiliary}
Let~$\delta$ be a real number ${0<\delta\le 1/2}$ and let~$N$ be a positive integer. 
There exists a non-zero polynomial ${G(X,Y) \in \QQ[X,Y]}$ satisfying
${\deg_X G \leq N}$, ${\deg_Y G \leq n-1}$,
\begin{align}
\label{eauxmult}
\nu_x (G (x, y(x))) &\geq (1 - \delta ) Nn,\\
\label{eauxh}
\height_p(G) &\leq \delta^{-1} nN(\height(\A) + 3) .
\end{align}
\end{proposition}

\paragraph{Proof}
It is quite analogous to the proof of Proposition~3.1 in~\cite{Ab08}. 
Condition~\eqref{eauxmult} is equivalent to a system of ${(1 - \delta ) Nn}$ linear equations in the ${n(N+1)}$ coefficients of~$G$. Each coefficient of each linear equation is a coefficient of~$x^k$, for ${k\le Nn}$, one of the series $y(x)^\ell$ for ${\ell=0,\ldots, n-1}$.

Using~\eqref{eeis1} and Lemma~\ref{lPS}, we estimate the height of every equation as ${nN\height(\A)+ (Nn+n)\log 2}$. Corollary~\ref{cdp} implies now that we can find a non-zero solution of our system of height at most 
$$
\delta^{-1}(nN\height(\A)+ (Nn+n)\log2)+ \frac12\delta^{-1}\log(nN).
$$
This is smaller than the right-hand side of~\eqref{eauxh}. \qed

\bigskip

\subsection{Upper Bound}
\label{ssup}

Now we can  obtain an upper bound for $\height_T(\alpha)$ in terms of ${\height(\alpha)}$. 

\begin{proposition}[Upper bound for $\height_T(\alpha)$]
\label{pupp}
Let~$\delta$ satisfy ${0<\delta \le 1/2}$. Then 
we have either
\begin{equation}
\label{eledelta}
\height(\alpha)\le 10\delta^{-2}mn^4(\height_p(F)+5),
\end{equation}
 or 
\begin{equation}
\label{eupp}
n\height_T(\alpha)\le (1+4\delta)\height(\alpha)+ 8\delta^{-1} n(\height(\A) + 10)+ \height_p(F).
\end{equation}
\end{proposition}

\paragraph{Proof}
Fix a positive integer~$N$, to be specified later,  and let $G(X,Y)$
be the auxiliary polynomial  introduced in Proposition~\ref{pAuxiliary}.   Extending the field~$\K$, we may  assume that ${G(X,Y)\in \K[X,Y]}$. 
We may also assume that~$G$ has a coefficient equal to~$1$; in particular, ${|G|_v\ge 1}$ for all ${v\in M_\K}$, where we denote by~$|G|_v$ the maximum of $v$-adic norms of coefficients of~$G$.

The series ${z(x)=G(x,y(x))\in \K[[x]]}$ can be written as
$$
z(x)=\sum_{k=\eta}^\infty b_kx^k
$$
with ${\eta\ge (1-\delta)Nn\ge\frac12Nn}$ (recall that ${\delta\le 1/2}$). Again using~\eqref{eeis1} and Lemma~\ref{lPS}, we estimate the coefficients~$b_k$ as follows: for ${v<\infty}$ we have ${|b_k|_v\le |G|_vA_v^k}$, and for ${v\mid\infty}$ we have 
${|b_k|_v\le n(N+1)2^{k+n-1}|G|_vA_v^k}$. Since for ${k\ge\eta\ge \frac12Nn}$ we have ${n(N+1)2^{k+n-1}\le 8^k}$, we obtain the estimate 
\begin{equation}
\label{ebkv}
|b_k|\le 
\begin{cases}
|G|_vA_v^k, & v<\infty,\\
|G|_v(8A_v)^k, & v\mid\infty.
\end{cases}
\qquad (v\in M_k, \quad, k\ge \eta). 
\end{equation}
Now we distinguish two cases.

\paragraph{Case~1: ${G(\alpha,\beta)=0}$}

In this case we have ${F(\alpha,\beta)=G(\alpha,\beta)=0}$.  We want to apply Lemma~\ref{lsys}; for this, we have to verify that  polynomials~$F$ and~$G$ do not have a common factor. This is indeed the case, because~$F$ is absolutely  irreducible, and ${\deg_YG<\deg_YF}$. 

Lemma~\ref{lsys}, combined with~\eqref{eauxh} and~\eqref{eeis2}, gives
\begin{align}
\height(\alpha)&\le  n \height_p (G) + (n-1) \height_p F + (m (n-1) + Nn) + (2n-1) \log (2n-1)+\log2\nonumber\\
&\le \delta^{-1}Nn^2(\height(\A)+6) + (n-1)(\height_p(F)+m)\nonumber\\
\label{eleal}
&\le 5\delta^{-1}Nn^3(\height_p(F)+5)+mn. 
\end{align}
Below, after specifying~$N$, we will see that this is sharper than~\eqref{eledelta}. 

\paragraph{Case~2: ${G(\alpha,\beta)=\gamma\ne0}$}
To treat this case it will be convenient to use, instead of the set~$T$,   a slightly smaller subset~$\tilT$,
 consisting of ${v\in T}$ satisfying 
\begin{equation*}
|\alpha|_v <
\begin{cases}
A_v^{-1}, &v<\infty,\\
(16A_v)^{-1}, &v\mid\infty. 
\end{cases}
\end{equation*}
We have clearly
\begin{equation}
\label{ettt}
0\le\height_T(\alpha)-\height_\tilT(\alpha)\le \height(\A)+\log 16,
\end{equation}
and~\eqref{ebkv} implies the estimate
\begin{equation}
\label{ebkak}
|b_k\alpha^k|_v<
\begin{cases}
|G|_vA_v^\eta|\alpha|_v^\eta, &v<\infty,\\
|G|_v (8A_v)^\eta |\alpha|_v^\eta \cdot (1/2)^{k-\eta}, &v\mid\infty. 
\end{cases}
\qquad (v\in \tilT,\quad k\ge\eta). 
\end{equation}

Recall that for ${v\in T}$, the series $y(x)$ converges $v$-adically to~$\beta$ at ${x=\alpha}$. Hence the same holds true for ${v\in \tilT}$. It follows that, for ${v\in \tilT}$, the series ${z(x)=G(x,y(x))}$ converges $v$-adically to\footnote{For archimedean~$v$ to make this conclusion  we need absolute convergence of $y(x)$ at ${x=\alpha}$, which is obvious for for ${v\in \tilT}$.}  ${G(\alpha,\beta)=\gamma}$. 

Using~\eqref{ebkak}, we can estimate $|\gamma|_v$ for ${v\in \tilT}$:
\begin{equation*}
|\gamma|_v<
\begin{cases}
|G|_vA_v^\eta|\alpha|_v^\eta, &v<\infty,\\
2|G|_v (8A_v)^\eta |\alpha|_v^\eta , &v\mid\infty. 
\end{cases}
\qquad (v\in \tilT,\quad k\ge\eta). 
\end{equation*}
Using this and remembering that ${|G|_v\ge 1}$ for all~$v$, we obtain the following lower estimate for ${\height(\gamma)}$:
\begin{align*}
\height(\gamma)&\ge\height_{\tilT}(\gamma)\\
&\ge \eta \height_\tilT(\alpha)-\height_p(G) -\eta\height(\A)-\eta\log16-\log2\\
&\ge Nn(1-\delta)\height_\tilT(\alpha) - 2\delta^{-1} nN(\height(\A) + 6).
\end{align*}
Combining this with~\eqref{ettt}, we obtain 
\begin{equation}
\label{elohg}
\height(\gamma)\ge Nn(1-\delta)\height_T(\alpha) - 3\delta^{-1} nN(\height(\A) + 6).
\end{equation}

On the other hand, using Lemma~\ref{lmaxG} it is easy to bound $\height(\gamma)$ from above. Indeed, part~\ref{iba} of this lemma implies that 
$$
\height(\beta) \leq \height_p(F) + m\height(\alpha) +n +\log (m+1), 
$$
and part~\ref{ifab} implies that 
$$
\height(\gamma) \le \height_a (G) + N\height(\alpha) +(n-1)\height(\beta) + \log((N+1)n).
$$
Since~$G$ has a coefficient equal to~$1$, we have ${\height_a(G)=\height_p(G)\le \delta^{-1} nN(\height(\A) + 3)}$. Hence 
\begin{align*}
\height(\gamma) &\le \height_p(G) + N\height(\alpha) +(n-1)(\height_p(F) + m\height(\alpha) +n +\log (m+1)) + \log((N+1)n)\\
&\le (N+mn)\height(\alpha)  + \delta^{-1} nN(\height(\A) + 4) +n\height_p(F)+n^2+n\log(m+1) .
\end{align*}
Combining this with~\eqref{elohg} and dividing by~$N$, we obtain
\begin{equation}
\label{emessy}
n(1-\delta)\height_T(\alpha)\le \left(1+\frac{mn}N\right)\height(\alpha)  + 4\delta^{-1} n(\height(\A) + 6) +N^{-1}(n\height_p(F)+n^2+n\log(m+1)). 
\end{equation}

\paragraph{Completing the proof of Proposition~\ref{pupp}}
Now it is the time to specify~$N$: we set ${N=\lceil\delta^{-1}mn\rceil}$. With this choice of~$N$, inequality~\eqref{eleal} is indeed sharper than~\eqref{eledelta}, and inequality~\eqref{emessy} implies the following:
$$
n(1-\delta) \height_T(\alpha) \le (1+\delta)\height(\alpha)+ 4\delta^{-1} n(\height(\A) + 10)+ \delta\height_p(F). 
$$
Since $\delta\le 1/2$, this is sharper than~\eqref{eupp}. \qed

\subsection{Lower Bound}
\label{sslow}
Our next objective is a lower bound for $\height_T(\alpha)$. We will see that it easily follows from the upper bound. 

\begin{proposition}[Lower bound for $\height_T(\alpha)$]
\label{plow}
Let~$\delta$ satisfy ${0<\delta \le 1/2}$. Then 
we have either~\eqref{eledelta} or 
\begin{equation}
\label{elow}
n\height_T(\alpha)\ge (1-4n\delta)\height(\alpha)- 9\delta^{-1} n^2(\height(\A) + 10)- n\height_p(F).
\end{equation}
\end{proposition}

\paragraph{Proof}
Remark first of all that we may assume that the polynomial $F(\alpha,Y)$ is of degree~$n$ and separable. Indeed, if this is not the case, then ${R_F(\alpha)=0}$, where
 $R_F(X)$ is the $Y$-resultant of $F(X,Y)$ and its $Y$-derivative $F'_Y(X,Y)$. In this case, the joint application of Lemmas~\ref{lSchmidt2} and~\ref{Yuri2} gives
\begin{equation*}
\height(\alpha)\le 2n\height_p (F) + 2n \log ((m+1)(n+1)\sqrt{n})+\log 2,
\end{equation*}
sharper than~\eqref{eledelta}. 

Thus, $F(\alpha,Y)$ has~$n$ distinct roots in~$\bar\Q$, one of which is~$\beta$;  we denote them ${\beta_1=\beta, \beta_2, \ldots, \beta_n}$. Extending the field~$\K$, we may assume that ${\beta_1, \ldots, \beta_n\in \K}$.

Set ${S=\{v\in M_\K: |\alpha|_v<1\}}$. 
For ${i=1,\ldots,n}$ we let~$T_i$ be the set of ${v\in S}$ such that $y(x)$ converges $v$-adically to~$\beta_i$ at ${x=\alpha}$; in particular, ${T_1=T}$. The sets ${T_1, \ldots, T_n}$ are clearly disjoint, and we have
\begin{equation}
\label{ests}
S\supset T_1\cup\ldots\cup T_n\supset \tilS,
\end{equation}
where~$\tilS$ consists of ${v\in S}$ for which ${|\alpha|_v<A_v^{-1}}$. The left inclusion in~\eqref{ests} is trivial, and to prove the right one just observes that for every ${v\in\tilS}$, the series $y(x)$ absolutely converges $v$-adically at ${x=\alpha}$, and, since ${F(x,y(x))=0}$, the sum must be a root of $F(\alpha,Y)$.

Clearly, 
$$
0\le\height(\alpha)-\height_\tilS(\alpha)=\height_S(\alpha)-\height_\tilS(\alpha)\le \height(\A). 
$$
It follows that
$$
\height_{T_1}(\alpha)+\cdots+\height_{T_n}(\alpha)\ge \height_\tilS(\alpha)\ge \height(\alpha)-\height(\A). 
$$
Now observe that the upper bound~\eqref{eupp} holds true with~$T$ replaced by any~$T_i$:
$$
n\height_{T_i}(\alpha)\le (1+4\delta)\height(\alpha)+ 8\delta^{-1} n(\height(\A) + 10)+ \height_p(F) \qquad (i=1, \ldots, n). 
$$
The last two inequalities imply that 
$$
n\height_{T}(\alpha)=n\height_{T_1}(\alpha)\ge n(\height(\alpha)-\height(\A))-(n-1)((1+4\delta)\height(\alpha)+ 8\delta^{-1} n(\height(\A) + 10)+ \height_p(F)),
$$
which easily transforms into~\eqref{elow}. \qed

\subsection{Proof of the ``Main Lemma''}
Using Propositions~\ref{pupp} and~\ref{plow} with ${\delta=\eps/4}$  and dividing by~$n$, we obtain that either~\eqref{eml1} holds, or 
$$
\left|\height_T(\alpha)-\frac{\height(\alpha)}n\right|\le \eps \height(\alpha)+ 40\eps^{-1} n(\height(\A) + 10)+ \height_p(F).
$$ 
Combining this with~\eqref{eeis2}, we obtain~\eqref{eml2}. \qed

\subsection{``Ramified Main Lemma''}
\label{ssrml}
We will actually need a slightly more general statement, allowing ramification in the series $y(x)$. The set-up is  as before, except that  now we consider the series
$$
y(x)=\sum_{k=1}^\infty a_kx^{k/e}\in \K[[x^{1/e}]]
$$
satisfying ${F(x,y(x))=0}$. We fix an $e$-th root $\alpha^{1/e}$ and we will assume that it belongs to~$\K$. We will now say that the series $y(x)$  converges $v$-adically to~$\beta$ at~$\alpha$ if the series $y(x^e)$ converges $v$-adically to~$\beta$ at $\alpha^{1/e}$. (Of course, this depends on the particular choice of the root $\alpha^{1/e}$.) We again define~$T$ as the set of all ${v\in S}$ for which $y(x)$  converges $v$-adically to~$\beta$ at~$\alpha$. 

\begin{lemma}[``Ramified Main Lemma'']
\label{lrml}
Let~$\eps$ satisfy ${0<\eps\le 1}$. Then  we have either
\begin{equation}
\label{erml1}
\height(\alpha)\le 200\eps^{-2}me^2n^4(\height_p(F)+5),
\end{equation}
or
\begin{equation}
\label{erml2}
\left|\frac{\height(\alpha)}n-\height_T(\alpha)\right|\le \eps\height(\alpha)+200\eps^{-1} en^2(\height_p(F)+2\log(mn) + 10). 
\end{equation}
\end{lemma}

\paragraph{Proof}
The proof is by reduction to the unramified case. Apply Lemma~\ref{lml} to the polynomial $F(X^e,Y)$, the series $y(x^e)$ and the number $\alpha^{1/e}$. We obtain that either
\begin{equation*}
\height(\alpha^{1/e})\le 200\eps^{-2}men^6(\height_p(F)+5),
\end{equation*}
or
\begin{equation*}
|\height(\alpha^{1/e})-n\height_T(\alpha^{1/e})|\le\eps \height(\alpha^{1/e})+200\eps^{-1} n^4(\height_p(F)+\log(men) + 10). 
\end{equation*}
These estimates easily transform into~\eqref{erml1} and~\eqref{erml2}, respectively, using that
$$
\height(\alpha^{1/e})=e^{-1}\height(\alpha), \quad \height_T(\alpha^{1/e})=e^{-1}\height_T(\alpha),\quad e\le n. \eqno\square
$$

\section{Proof of the Main Theorem}
\label{sproof}

In this section we prove Theorem~\ref{MainResult}. First of all, we investigate the relation between $\height_T(\alpha)$ and $\lgcd_T(\alpha,\beta)$, where~$T$ is defined as in Section~\ref{sspr}.  

\subsection{Comparing \texorpdfstring{$\height_T(\alpha)$}{hT(a)} and \texorpdfstring{$\lgcd_T(\alpha,\beta)$}{lgcdT(a,b)}}
In this subsection we retain the set-up of Subsection~\ref{ssml}, except that we allow ramification in the series $y(x)$, as we did in Subsection~\ref{ssrml}. Thus, in this subsection:
\begin{itemize}
\item
$\K$ is a number field;
\item
$F(X,Y)\in\K[X,Y]$ is an absolutely irreducible polynomial;
\item
$\alpha,\beta\in \K$ satisfy ${F(\alpha,\beta)=0}$;
\item
${y(x)=\sum_{k=1}^\infty a_kx^{k/e}\in \K[[x^{1/e}]]}$
satisfies ${F(x,y(x))=0}$;

\item
${T\subset M_\K}$ is the set of all ${v\in M_\K}$ such that ${|\alpha|_v<1}$ and $y(x)$ converges $v$-adically at~$\alpha$ to~$\beta$.
\end{itemize}
The $v$-adic convergence is understood in the same sense as in Subsection~\ref{ssrml}: we fix an $e$-th root $\alpha^{1/e}$, assume that it belongs to~$\K$ and  and define $v$-adic convergence of $y(x)$ to~$\beta$ at~$\alpha$ as $v$-adic convergence of $y(x^e)$ to~$\beta$ at~$\alpha^{1/e}$.

Let~$\kappa$ be the smallest~$k$ such that ${a_k\ne 0}$; by the assumption, ${\kappa>0}$. Then we have ${\nu_x(y)=\kappa/e}$ and
$$
y(x)=\sum_{k=\kappa}^\infty a_kx^{k/e}
$$
with ${a_\kappa\ne 0}$. In this subsection we prove that $\lgcd_T(\alpha,\beta)$ can be approximated by $\min\{1,\kappa/e\}\height_T(\alpha)$.

\begin{proposition}
\label{pnew}
In the above set-up we have 
\begin{equation}
\label{edhgcd}
\left|\lgcd_T(\alpha,\beta)-\min\{\kappa/e,1\}\height_T(\alpha)\right|\le 30n\kappa\height_p (F) + 30n\kappa \log (nm) + 15en. 
\end{equation}
\end{proposition}

This statement corresponds to Proposition~3.6 in~\cite{Ab08}. Our proof is, however, much more involved, in particular  because Abouza\"id did not need the lower estimate.

\paragraph{Proof}
Let ${\A = \left( A_v \right)_{v \in M_{\K}}}$ be the $M_\K$-divisor from Corollary~\ref{cBB}. For the reader's convenience, we reproduce here~\eqref{eeis1} and~\eqref{eeis2}:
\begin{align*}
|a_k|_v &\leq A_v^{k/e} \qquad(v\in M_\K, \quad k\ge 1),\\
\height (\A) &\leq 4n\height_p (F) + 3n \log (nm) + 10en.
\end{align*}
As we already did several times in Section~\ref{sspr}, it will be convenient to replace~$T$ by a smaller subset. Thus, let~$\tilT$ consist of ${v\in T}$ satisfying
\begin{equation}
\label{eleka}
|\alpha|_v< 
\begin{cases}
A_v^{-\kappa-1}\min\{1,|a_\kappa|_v\}^e, & v<\infty,\\
(1/4)^eA_v^{-\kappa-1}\min\{1,|a_\kappa|_v\}^e, & v<\infty.
\end{cases}
\end{equation}
(Attention: this is not the same~$\tilT$ as in Subsection~\ref{ssup}!) 
Clearly,
$$
0\le\height_T(\alpha)-\height_\tilT(\alpha)\le (\kappa+1)\height(\A)+e\height_{T\smallsetminus\tilT}(a_\kappa).
$$
Using~\eqref{eeis1} we estimate ${\height(a_\kappa)\le (\kappa/e)\height(\A)}$. We obtain
\begin{equation}
\label{edifhe}
0\le\height_T(\alpha)-\height_\tilT(\alpha)\le (\kappa+1)\height(\A)\le3\kappa\height(\A)+e\log 4,
\end{equation}
where for the latter estimate we use ${\kappa\ge 1}$. 
In particular, 
\begin{equation}
\label{edifgcd}
0\le\lgcd_T(\alpha,\beta)-\lgcd_\tilT(\alpha,\beta)\le 3\kappa\height(\A)+e\log4.
\end{equation}

After this preparation, we can now proceed with the proof. For every ${v\in \tilT}$ we want to obtain an estimate of the form 
${c_v|\alpha|_v^{\kappa/e}\le|\beta|_v\le c_v'|\alpha|_v^{\kappa/e}}$, 
where $c_v$ and $c'_v$ are some quantities not depending on~$\alpha$.   

\paragraph{Upper estimate for $|\beta|_v$.}
This is easy. It follows from~\eqref{eleka} that
$$
|\alpha|_v< 
\begin{cases}
A_v^{-1}, & v<\infty,\\
(4^eA_v)^{-1}, & v<\infty.
\end{cases}
$$
From this and~\eqref{eeis1} we deduce that 
\begin{equation}
\label{eakake}
|a_k\alpha^{k/e}|_v<
\begin{cases}
A_v^{\kappa/e}|\alpha|_v^{\kappa/e},& v<\infty,\\
A_v^{\kappa/e}|\alpha|_v^{\kappa/e}\cdot (1/4)^{k-\kappa},& v\mid\infty
\end{cases}
\qquad (k\ge \kappa).
\end{equation}
Hence
$$
|\beta|_v<
\begin{cases}
A_v^{\kappa/e}|\alpha|_v^{\kappa/e},& v<\infty,\\
2A_v^{\kappa/e}|\alpha|_v^{\kappa/e},& v\mid\infty. 
\end{cases}
$$

\paragraph{Lower estimate for $|\beta|_v$.}
The lower estimate is slightly more subtle. First, we bound the difference ${\beta-a_\kappa\alpha^{\kappa/e}}$ from above using~\eqref{eleka}. 

Similarly to~\eqref{eakake}, we have 
$$
|a_k\alpha^{k/e}|_v<
\begin{cases}
A_v^{(\kappa+1)/e}|\alpha|_v^{(\kappa+1)/e},& v<\infty,\\
A_v^{(\kappa+1)/e}|\alpha|_v^{(\kappa+1)/e}\cdot (1/4)^{(k-\kappa-1)/e},& v\mid\infty
\end{cases}
\qquad (k\ge \kappa+1).
$$
Hence, presenting ${\beta-a_\kappa\alpha^{\kappa/e}}$ as the $v$-adic sum of the series
$$
y(x)-a_\kappa x^{\kappa/e}=\sum_{k=\kappa+1}^\infty a_kx^{k/e}
$$
at ${x=\alpha}$, we obtain the estimate
$$
|\beta-a_\kappa\alpha^{\kappa/e}|_v< 
\begin{cases}
A_v^{(\kappa+1)/e}|\alpha|_v^{(\kappa+1)/e},& v<\infty,\\
2A_v^{(\kappa+1)/e}|\alpha|_v^{(\kappa+1)/e},& v\mid\infty.
\end{cases}
$$
Combining this with~\eqref{eleka}, we find
$$
|\beta-a_\kappa\alpha^{\kappa/e}|_v< 
\begin{cases}
\min\{|a_\kappa|_v,1\}|\alpha|_v^{\kappa/e},& v<\infty,\\
(1/2)\min\{|a_\kappa|_v,1\}|\alpha|_v^{\kappa/e},& v\mid\infty.
\end{cases}
$$
Hence
$$
|\beta|_v\ge 
\begin{cases}
\min\{|a_\kappa|_v,1\}|\alpha|_v^{\kappa/e},& v<\infty,\\
(1/2)\min\{|a_\kappa|_v,1\}|\alpha|_v^{\kappa/e},& v\mid\infty,
\end{cases}
$$
the lower estimate we were seeking.

\paragraph{Completing the proof of Proposition~\ref{pnew}}
Thus, we proved that
\begin{equation}
\label{edouble}
c_v|\alpha|_v^{\kappa/e}\le|\beta|_v\le c_v'|\alpha|_v^{\kappa/e}, \\
\end{equation}
with 
\begin{equation*}
c_v=
\begin{cases}
\min\{|a_\kappa|_v,1\},& v<\infty,\\
(1/2)\min\{|a_\kappa|_v,1\},& v\mid\infty,
\end{cases}, \qquad
c'_v=
\begin{cases}
A_v^{\kappa/e},& v<\infty,\\
2A_v^{\kappa/e},& v\mid\infty. 
\end{cases}
\end{equation*}
From~\eqref{edouble} we deduce that for ${v\in \tilT}$
$$
c_v|\alpha|_v^{\min\{\kappa/e,1\}}\max\{|\alpha|_v,|\beta|_v\}\le c_v'|\alpha|_v^{\min\{\kappa/e,1\}}.
$$
(We use here the obvious inequality ${c_v\le1\le c'_v}$.) Hence
$$
-(\kappa/e)\height(\A)-\log2\le \lgcd_\tilT(\alpha,\beta)-\min\{\kappa/e,1\}\height_\tilT(\alpha)\le \height(a_\kappa)+\log2.
$$
Since ${\height(a_\kappa)\le (\kappa/e)\height(\A)}$, this implies 
$$
|\lgcd_\tilT(\alpha,\beta)-\min\{\kappa/e,1\}\height_\tilT(\alpha)|\le (\kappa/e)\height(\A)+\log2,
$$
which, together with~\eqref{edifhe} and~\eqref{edifgcd} gives
$$
|\lgcd_\tilT(\alpha,\beta)-\min\{\kappa/e,1\}\height_\tilT(\alpha)|\le 7\kappa\height(\A)+4e.
$$
Combining this with~\eqref{eeis1}, we obtain~\eqref{edhgcd}. \qed

\subsection{Proving Theorem~\ref{MainResult}}
Now we are fully equipped for the proof of our main result. We want to show that, assuming
\begin{equation}
\label{eassum}
\height(\alpha)\ge 200\eps^{-2}mn^6(\height_p(F)+5),
\end{equation}
we have 
\begin{equation}
\label{eresu}
\left| \frac{\lgcd (\alpha, \beta)}{r} - \frac{\height(\alpha)}{n} \right| \leq
\frac{1}{r} \left(   \eps\height(\alpha)+4000\eps^{-1} n^4(\height_p(F)+\log(mn) + 1) +30n^2m(\height_p (F) + \log (nm))  \right).
\end{equation}

Write 
${F(X,Y)=f_n(X)Y^n+\cdots+f_0(X)}$.
According to Proposition~\ref{ppe} we have
\begin{equation*}
F(x,Y)=f_n(x)(Y-y_1)\cdots(Y-y_n).
\end{equation*}
where
$$
y_i=\sum_{k=\kappa_i}^\infty a_{ik}x^{k/e_i}\in \K((x^{1/e_i}))  \qquad (i=1, \ldots, n). 
$$
We assume that  ${a_{i\kappa_i}\ne 0}$ for ${i=1, \ldots, n}$, so that ${\kappa_i/e_i=\nu_x(y_i)}$. 

Denoting by~$\ell$ the number of indexes~$i$ such that ${\kappa_i>0}$, we may  assume that ${\kappa_1, \ldots, \kappa_\ell>0}$ and ${\kappa_{\ell+1}, \ldots, \kappa_n\le 0}$. Propositions~\ref{prke}  implies that 
\begin{equation}
\label{er=}
r=\sum_{i=1}^\ell\min\{1,\kappa_i/e_i\},
\end{equation}
and Proposition~\ref{pfs} implies that ${f_\ell(0)\ne 0}$. We may normalize polynomial ${F(X,Y)}$ to have
\begin{equation*}
f_\ell(0)=1. 
\end{equation*}
In particular, ${|F|_v\ge 1}$ for every ${v\in M_\K}$, where $|F|_v$ denotes the maximum of $v$-adic norms of the coefficients of~$F$, and also ${\height_p(F)=\height_a(F)}$.

Set ${E=\lcm(e_1,\ldots,e_\ell )}$ and fix an $E$-th root $\alpha^{1/E}$. This fixes uniquely the roots ${\alpha^{1/e_1}, \ldots, \alpha^{1/e_\ell}}$. 
Extending the field~$\K$ we may assume that the coefficients of the series ${y_1, \ldots, y_\ell}$  belong to~$\K$, and the same is true for $\alpha^{1/E}$ (and hence for  ${\alpha^{1/e_1}, \ldots, \alpha^{1/e_\ell}}$ as well). Having fixed the root ${\alpha^{1/e_i}\in \K}$, we may define $v$-adic convergence of $y_i$ at~$\alpha$, see Subsection~\ref{ssrml}. 

Extending further the field~$\K$, we may assume that it contains all the roots of the polynomial $F(\alpha,Y)$. Hence, if one of the series ${y_1,\ldots,y_\ell}$ converges $v$-adically at~$\alpha$ (and if the convergence is absolute in the archimedean case), then the sum must belong to~$\K$.

Consider the following subsets of~$M_\K$:
\begin{align*}
S&=\{v\in M_K: |\alpha|_v<1\},\\
T_i&=\{v\in S: \text{the series $y_i$ converges $v$-adically to~$\beta$ at~$\alpha$}\} \qquad (i=1, \ldots, \ell ). 
\end{align*}
(These sets are not the same~$T_i$ as in Subsection~\ref{sslow}!)  

We have clearly ${\lgcd(\alpha,\beta)=\lgcd_S(\alpha,\beta)}$. If we manage to show that 
the sets~$T_i$ are pairwise disjoint, and that
$\height_{S\smallsetminus(T_1\cup\cdots\cup T_\ell)}(\beta)$ is ``negligible'', 
then joint application of Lemma~\ref{lrml}, Proposition~\ref{pnew} and identity~\eqref{er=} would prove Theorem~\ref{MainResult}. We will argue like this, only with the sets~$T_i$  replaced by slightly smaller subsets. 

Let  ${\A_i=(A_{iv})_{v\in M_\K}}$ be the $M_\K$-divisor for the series~$y_i$ given by Corollary~\ref{cBB}. Define the $M_\K$-divisor ${\A=(A_{v})_{v\in M_\K}}$ by
$$
A_v=\max\{A_{1v},\ldots,A_{\ell v}\}\qquad (v\in M_\K). 
$$
We have clearly 
\begin{align}
\label{eakbo} |a_{ki}|_v &\leq A_v^{k/e} \qquad(v\in M_\K, \quad 1\le i\le \ell ,  \quad k\ge \kappa_i),\\
\height (\A) &\leq \height(\A_1)+\cdots+\height(\A_\ell)\nonumber\\
&\le  4n^2\height_p (F) + 3n^2 \log (nm) + 10n^3. \nonumber
\end{align}
Now let~$\tilS$ consist of ${v\in S}$ satisfying
\begin{equation}
\label{ealbo}
|\alpha|_v<
\begin{cases}
|F|_v^{-n}A_v^{-1}, & v<\infty,\\
((n+1)2^{n+3}|F|_v)^{-n}A_v^{-1}, & v\mid\infty, 
\end{cases}
\end{equation}
and set ${\tilT_i=T_i\cap\tilS}$. (This is not the same~$\tilS$ that in Subsection~\ref{sslow}!) Clearly,
\begin{align}
0\le \lgcd(\alpha,\beta)-\lgcd_\tilS(\alpha,\beta)&\le \height(\alpha)-\height_\tilS(\alpha)\nonumber\\ &=\height_{S\smallsetminus\tilS}(\alpha)\nonumber\\ &\le \height(\A)+n\height_p(F)+\log ((n+1)2^{n+3})\nonumber\\
\label{esisi}
&\le  5n^2\height_p (F) + 3n^2 \log (nm) + 15n^3,\\
0\le \lgcd_{T_i\smallsetminus\tilT_i}(\alpha,\beta) &\le\height_{S\smallsetminus\tilS}(\alpha)\nonumber\\ 
\label{etiti}
&\le  5n^2\height_p (F) + 3n^2 \log (nm) + 15n^3 \qquad (i=1, \ldots, \ell). 
\end{align}
Here we used the equality ${\height_p(F)=\height_a(F)}$. 

Mention also that for ${v\in \tilS}$, we have ${|\alpha|_v<A_v^{-1}}$, which implies that the series ${y_1, \ldots, y_\ell}$ converge $v$-adically at~$\alpha$ in the completion~$\K_v$,  the convergence being absolute when~$v$ is archimedean. Hence, as we have seen above, the sum must belong to~$\K$.

\begin{proposition}
\label{pdisj}
The sets ${\tilT_1, \ldots, \tilT_\ell}$ pairwise disjoint. Furthermore, if ${v\in \tilS}$ but ${v\notin \tilT_1\cup\ldots\cup \tilT_\ell}$ then 
\begin{equation}
\label{ebge}
|\beta|_v\ge 
\begin{cases}
|F|_v^{-1},&v<\infty,\\
((n+1)2^{n+2}|F|_v)^{-1},&v\mid\infty.
\end{cases}
\end{equation}

\end{proposition}

\paragraph{Proof}
The polynomial 
$$
Q(Y)=(Y-y_1)\cdots(Y-y_\ell)\in \K[[x^{1/E}]][Y]. 
$$
divides $F(x,Y)$ in the ring $\K((x^{1/E}))[Y]$. By Gauss' Lemma, $Q(Y)$ divides $F(x,Y)$ in the ring  $\K[[x^{1/E}]][Y]$ as well. Moreover, writing ${F(x,Y)=Q(Y)U(Y)}$ with
$$
U(Y)= f_n(x)Y^{n-\ell }+u_{n-\ell -1}Y^{n-\ell -1}+\cdots+u_0\in \K[[x^{1/E}]](Y),
$$
the coefficients ${u_0,\ldots, u_{n-\ell -1}}$  belong to the ring\footnote{This is a consequence of the general algebraic property: let~$R$ be a commutative ring,~$R'$ a subring and ${Q(Y),F(Y)\in R'[Y]}$, the polynomial~$Q$ being monic; assume that ${Q\mid F}$ in $R[Y]$; then ${Q\mid F}$ in $R'[Y]$. Indeed, denoting by~$a$ the leading coefficient of~$F$, the polynomial~$Q$ divides ${G=F-aY^{\deg F-\deg Q}Q}$ in $R[Y]$, and ${\deg G<\deg F}$, so by induction ${Q\mid G}$ in $R'[Y]$.} $\K[x,y_1,\ldots, y_\ell]$. Recall that for ${v\in \tilS}$ the series ${y_1, \ldots, y_\ell}$ converge $v$-adically at~$\alpha$ in the field~$\K$,  the convergence being absolute when~$v$ is archimedean.  Hence so  do the coefficients of~$U$. 

Fix ${v\in \tilS}$ and write
$$
F(\alpha,Y)=(Y-y_1(\alpha))\cdots(Y-y_\ell(\alpha))(f_n(\alpha)Y^{n-\ell }+u_{n-\ell -1}(\alpha)Y^{n-\ell -1}+\cdots+u_0(\alpha)),
$$
where  ${y_1(\alpha),\ldots,y_\ell(\alpha)\in\K}$ the $v$-adic sum of the corresponding series at~$\alpha$, and similarly for  ${u_{n-\ell -1}(\alpha),\ldots,u_0(\alpha)}$. We claim that $F(\alpha,Y)$ is a separable polynomial of degree~$n$; indeed, if this is not the case, then, as we have seen in Subsection~\ref{sslow}, our~$\alpha$ must satisfy~\eqref{eresu}, which contradicts~\eqref{eassum}. 

Now if ${v\in T_i\cap T_j}$ for ${i\ne j}$ then ${\beta=y_i(\alpha)=y_j(\alpha)}$, and $F(\alpha,Y)$ must have~$\beta$ as a double root, a contradiction.  This proves disjointedness of the sets~$\tilT_i$.  

Now assume that ${v\in \tilS}$ but ${v\notin \tilT_1\cup\ldots\cup \tilT_\ell}$. Then none of the sums ${y_1(\alpha),\ldots,y_\ell(\alpha)}$ is equal to~$\beta$; in other words ${y_1(\alpha),\ldots,y_\ell(\alpha),\beta}$ are ${\ell+1}$ distinct roots of the polynomial 
$$
P(Y)=F(\alpha, Y)=f_n(\alpha)Y^n+\cdots+f_0(\alpha). 
$$
We are going to use Lemma~\ref{lr+1}. Since ${f_\ell(0)=1}$ and 
$$
|\alpha|_v<
\begin{cases}
|F_v|^{-1},&v<\infty,\\
(2|F|_v)^{-1},&v\mid\infty,
\end{cases}
$$
we have 
$$
|f_\ell(\alpha)|_v\ge 
\begin{cases}
1,&v<\infty,\\
1/2,&v\mid\infty,
\end{cases},\qquad
|P|_v\le 
\begin{cases}
|F|_v,&v<\infty,\\
2|F|_v,&v\mid\infty.
\end{cases}
$$
Now Lemma~\ref{lr+1} implies that 
\begin{equation}
\label{emaxle}
\max\{|y_1(\alpha)|_v,\ldots,|y_\ell(\alpha)|_v,|\beta|_v\}\ge 
\begin{cases}
|F|_v^{-1},&v<\infty,\\
((n+1)2^{n+2}|F|_v)^{-1},&v\mid\infty.
\end{cases}
\end{equation}

On the other hand, we may estimate ${|y_i(\alpha)|_v}$ from above using~\eqref{eakbo} and~\eqref{ealbo}. In what follows we repeatedly use the inequality ${e_i\le n}$. Since
$$
|\alpha|_v<
\begin{cases}
A_v^{-1},&v<\infty,\\
(2^{e_i}A_v)^{-1},&v\mid\infty
\end{cases}  \qquad (i=1, \ldots,\ell),
$$
we have
$$
|a_k\alpha^{k/e_i}|_v<
\begin{cases}
(A_v|\alpha|_v)^{1/e_i},&v<\infty,\\
(A_v|\alpha|_v)^{1/e_i}\cdot (1/2)^{k-1},&v\mid\infty 
\end{cases} \qquad (k\ge 1,\quad i=1, \ldots,\ell),
$$
which implies
$$
|y_i(\alpha)|_v<
\begin{cases}
(A_v|\alpha|_v)^{1/e_i},&v<\infty,\\
2(A_v|\alpha|_v)^{1/e_i},&v\mid\infty
\end{cases} \qquad (i=1, \ldots,\ell).
$$
Now since  
$$
|\alpha|_v<
\begin{cases}
|F|_v^{-e_i}A_v^{-1}, & v<\infty,\\
((n+1)2^{n+3}|F|_v)^{-e_i}A_v^{-1}, & v\mid\infty
\end{cases}  \qquad (i=1, \ldots,\ell),
$$
we obtain finally
$$
|y_i(\alpha)|_v<
\begin{cases}
|F|_v^{-1},&v<\infty,\\
((n+1)2^{n+2}|F|_v)^{-1},&v\mid\infty
\end{cases}  \qquad (i=1, \ldots,\ell).
$$
Compared with~\eqref{emaxle}, this implies~\eqref{ebge}. The proposition is proved.\qed

\bigskip

An immediate consequence of the second statement of Proposition~\ref{pdisj} is the estimate 
\begin{equation}
\label{et''}
\lgcd_{\tilS\smallsetminus(\tilT_1\cup\ldots\cup\tilT_\ell)}\le\height_{\tilS\smallsetminus(\tilT_1\cup\ldots\cup\tilT_\ell)}(\beta)\le \height_p(F)+\log((n+1)2^{n+2})
\end{equation}
(we again use ${\height_a(F)=\height_p(F)}$).

\bigskip

Now we collect everything together to prove Theorem~\ref{MainResult}. According to Lemma~\ref{lrml}, condition~\eqref{eassum} implies that
$$
\left|\frac{\height(\alpha)}n-\height_{T_i}(\alpha)\right|\le \eps\height(\alpha)+200\eps^{-1} n^3(\height_p(F)+2\log(mn) + 10) \qquad (i=1, \ldots,\ell).
$$
Combining this with Proposition~\ref{pnew} and estimate~\eqref{etiti}, we obtain
\begin{align*}
\left|\min\bigl\{\frac{\kappa_i}{e_i},1\bigr\}\frac{\height(\alpha)}n-\lgcd_{\tilT_i}(\alpha,\beta)\right|&\le \eps\height(\alpha)+3000\eps^{-1} n^3(\height_p(F)+\log(mn) + 1)\\& +30nm\height_p (F) + 30nm \log (nm).  \qquad (i=1, \ldots,\ell).
\end{align*}
Summing up, using~\eqref{er=} and the disjointedness of the sets~$\tilT_i$, we obtain
\begin{align*}
\left|r\frac{\height(\alpha)}n-\lgcd_{\tilT_1\cup\ldots\cup\tilT_\ell}(\alpha,\beta)\right|&\le \eps\height(\alpha)+3000\eps^{-1} n^4(\height_p(F)+\log(mn) + 1)\\& +30n^2m\height_p (F) + 30n^2m \log (nm).  
\end{align*}
Finally, combining this with~\eqref{esisi} and~\eqref{et''}, we obtain~\eqref{eresu}. \qed

\end{document}